\documentclass[11pt,a4paper]{article}
\usepackage[T1]{fontenc}
\usepackage{lmodern}
\usepackage[margin=27mm]{geometry}
\usepackage{amsmath,amssymb,amsthm}
\usepackage{booktabs}
\usepackage{microtype}
\usepackage[hidelinks]{hyperref}
\usepackage{orcidlink}

\newtheorem{theorem}{Theorem}
\newtheorem{corollary}[theorem]{Corollary}
\newtheorem{lemma}[theorem]{Lemma}
\DeclareMathOperator{\mex}{mex}
\newcommand{\tail}{\operatorname{tail}_{10}}
\newcommand{\digits}[1]{\texttt{#1}}
\newcommand{\slc}[3]{#1[#2:#3]}
\newcommand{\rowcase}[4]{\(#1\) & \(#2\) & \(#3\) & \(#4\)\\}
\newcommand{\subtr}{\operatorname{SUB}}
\hypersetup{
  pdftitle={Quadratic periods and preperiods in subtraction games with one varying move},
  pdfauthor={Ruiqi Zhao, Qi Zhang, Zhihao Shi}
}
\title{\bfseries Quadratic periods and preperiods in subtraction games
       with one varying move}
\author{Ruiqi Zhao\textsuperscript{1,\textdagger\,\orcidlink{0009-0003-7217-1319}}\qquad
  Qi Zhang\textsuperscript{2}\qquad
  Zhihao Shi\textsuperscript{3}\\[4pt]
  \small\textsuperscript{1}The University of Tokyo, Tokyo, Japan\quad
  \href{mailto:richard-zhao@g.ecc.u-tokyo.ac.jp}
    {\texttt{richard-zhao@g.ecc.u-tokyo.ac.jp}}\\
  \small\textsuperscript{2}Beijing University of Posts and Telecommunications, Beijing, China\quad
  \href{mailto:zhangqi9008@bupt.edu.cn}
    {\texttt{zhangqi9008@bupt.edu.cn}}\\
  \small\textsuperscript{3}Anhui Science and Technology University, Anhui, China\quad
  \href{mailto:zhihao-shi@ahstu.edu.cn}
    {\texttt{zhihao-shi@ahstu.edu.cn}}\\[3pt]
  \small\textsuperscript{\textdagger}Corresponding author}
\date{}
\makeatletter
\renewcommand{\@maketitle}{%
  \begin{center}
    {\LARGE\@title\par}
    \vskip 9pt
    {\large\begin{tabular}[t]{@{}c@{}}\@author\end{tabular}\par}
  \end{center}}
\makeatother

\begin{document}
\maketitle
\vspace{-8pt}

\begin{abstract}
We give two infinite families of counterexamples to Zhang's conjectures on the period and
preperiod of the Sprague--Grundy sequence when a single move \(c\) is
adjoined to a fixed subtraction set \(S\). For \(S=\{4,5,10\}\) and
\(c=15t+13\) with \(t\ge2\), the least period is \(c+19\) and the
least preperiod is \((c^2+21c+8)/15\). For \(S=\{3,7,10\}\) and
\(c=90t+9\) with \(t\ge1\), the least period is
\((c^2+26c+615)/30\) and the least preperiod is \(19\) less than the
period. The latter quadratic period also occurs in the outcome
sequence, which is periodic from its first term. The proofs construct
the sequences from fixed words and reduce the mex recurrence to a
finite table of local identities.
\end{abstract}


\section{Introduction}

In the normal-play subtraction game \(\subtr(T)\), a move from \(n\)
subtracts an element \(s\) of a nonempty finite set \(T\) of positive integers,
provided \(s\le n\). The Sprague--Grundy (SG) sequence \(g_T\) is
determined by
\begin{equation}\label{eq:sg}
 g_T(n)=\mex\{g_T(n-s):s\in T,\ s\le n\},\qquad n\ge0,
\end{equation}
where \(\mex\) is the least nonnegative integer outside the displayed set.
The outcome sequence \(w_T(n)=\mathbf 1_{\{g_T(n)=0\}}\) marks losing
positions. This indicator is the complement of Zhang's
\cite[Remark 2]{zhang}; complementation preserves periods and preperiods.
For a sequence \(u\), let \(p(u)\) be its least positive eventual
period and let \(\ell(u)\) be the least \(N\ge0\) for which \(u\) is
periodic on \(n\ge N\). These quantities exist for \(g_T\) and \(w_T\):
after the first \(\max T\) terms, a window of \(\max T\) values of
\(g_T\) determines the next term, and each value is at most \(|T|\).

Zhang~\cite[Conjectures 1.2 and 1.3]{zhang} asked how \(p(g_T)\) and
\(\ell(g_T)\) depend on \(c\) when \(T=S\cup\{c\}\) and \(S\) is
fixed. Conjecture~1.2 predicts \(O(c)\) growth for both quantities.
Conjecture~1.3 predicts that both are eventually affine in \(c\)
on each residue class of some fixed modulus. Zhang also raised the
corresponding question for outcomes.
Earlier investigations of long periods include the families of
Alth\"ofer and B\"ultermann~\cite{althofer}, in which several moves
vary, and the seeded outcome sequences of Mikl\'os and Post~\cite{miklos}.
See Larsson and Saha~\cite{larsson} for a survey.

\begin{theorem}\label{thm:main}
Let \(T=S\cup\{c\}\). For integer \(t\) in the indicated ranges, the
least eventual periods and preperiods are as follows.
\begin{center}
\small
\renewcommand{\arraystretch}{1.24}
\begin{tabular}{@{}lcc@{}}
\toprule
 & \(S=\{4,5,10\}\) & \(S=\{3,7,10\}\)\\
\midrule
Adjoined move \(c\) & \(15t+13,\ t\ge2\) & \(90t+9,\ t\ge1\)\\
\(p(g_T)\) & \(c+19\) & \((c^2+26c+615)/30\)\\
\(\ell(g_T)\) & \((c^2+21c+8)/15\) & \(p(g_T)-19\)\\
\(p(w_T)\) & \(p(g_T)\) & \(p(g_T)\)\\
\(\ell(w_T)\) & \(c-1\) & \(0\)\\
\bottomrule
\end{tabular}
\end{center}
\end{theorem}

\begin{corollary}
For \(S=\{4,5,10\}\), the function \(\ell(g_{S\cup\{c\}})\)
is superlinear along \(c\equiv13\pmod{15}\).
For \(S=\{3,7,10\}\), the functions \(p(g_{S\cup\{c\}})\),
\(\ell(g_{S\cup\{c\}})\), and \(p(w_{S\cup\{c\}})\) are
superlinear along \(c\equiv9\pmod{90}\).
Each fixed set disproves both of Zhang's conjectures for SG sequences,
and the second disproves the outcome-period versions.
\end{corollary}

\begin{proof}
The indicated quantities have positive quadratic leading terms on
the stated progressions. A function eventually affine on each of
finitely many residue classes is \(O(c)\), since the slopes and
intercepts range over finite sets.
\end{proof}

\section{Local verification and minimality}\label{sec:method}

We index word positions from zero. Concatenation is written by
juxtaposition, and \(\slc{U}{a}{b}\) takes the letters at positions
\(a,\ldots,b-1\) of \(U\). We write \(|U|\) for its length,
\(U^r\) for \(r\)-fold concatenation, and \(U^\infty\) for infinite
repetition. The notation \(\tail(U)\) denotes the last ten letters of \(U\). Let \(G\) be a candidate for \(g_T\), specified
as a succession of rows. If a row \(B\) begins at \(b\ge10\), write
\(A\) for the ten letters immediately before the row. For the
positions \(n=b+z\ge c\) in this row, let \(F[z]=G(n-c)\).
In both families \(\max S=10<c\), so the mex condition is
\begin{equation}\label{eq:local}
 B[z]=\mex D_z,\qquad
 D_z=\{(AB)[10+z-s]:s\in S\}\cup\{F[z]\}.
\end{equation}
For \(z\) with \(b+z<c\), the corresponding entries of \(F\) are unused.
In both constructions \(G\) agrees with a periodic base word \(H^\infty\)
for \(n<c\). To check \(H^\infty\) for \(S\), it suffices to test
the legal options in \eqref{eq:sg} for \(0\le n<10+|H|\):
for \(n\ge10\), all moves in \(S\) are legal and both sides depend only
on \(n\bmod |H|\).
If these base checks and \eqref{eq:local} hold, induction on \(n\)
gives \(G=g_T\).

The tables below specify \(U,B,F\) with \(A=\tail(U)\), and \(z_0\)
denotes the first column of the row for which \(n\ge c\). Every table
entry satisfies \eqref{eq:local}. The following finite test spells out
the arithmetic. Its input words are tuples of nonnegative integers.
\begin{verbatim}
def row_ok(S, U, B, F, z0=0):
    AB = U[-10:] + B
    for z in range(z0, len(B)):
        options = {AB[10 + z - s] for s in S} | {F[z]}
        v = B[z]
        if v in options or any(i not in options for i in range(v)):
            return False
    return True
\end{verbatim}
The test uses the two defining properties of the mex. As an example,
for \(S=\{4,5,10\}\) in the \(U=B=H,\ F=\Phi(H,H)\) row of
Table~\ref{tab:first}, column \(8\) has option values \(1,0,1,3\)
and value \(2\). The checks in the two tables involve only fixed words
and the displayed two-valued choices. The parameter coverage is
established below by tracing the row and column of each option.

\begin{lemma}\label{lem:min}
Suppose \(P\) is a period of \(u\) from \(N\). If the counts of
specified cyclic factors in \(u[N:N+P]\), together with \(P\), have
greatest common divisor \(1\), then \(p(u)=P\).

If \(N>0\) and \(u(N-1)\ne u(N-1+P)\), then \(\ell(u)=N\).
\end{lemma}

\begin{proof}
Any eventual period \(q\) induces invariance of the cyclic word under
shift by \(q\), by applying the period relation on sufficiently late
copies of the word. Each orbit of starting positions under this
shift has size \(P/\gcd(P,q)\). The starting positions of each
specified cyclic factor form a union of these orbits, so its count
is divisible by this size. The gcd hypothesis forces
\(P/\gcd(P,q)=1\), hence \(P\mid q\) and \(p(u)=P\).
For the second assertion, suppose a period \(q\)
held from some \(M<N\). For \(n=N-1\), choose \(j\ge0\) large enough
that \(n+jq\ge N\). Then
\(u(n)=u(n+jq)=u(n+P+jq)=u(n+P)\), a contradiction.
\end{proof}

\section{The base set \texorpdfstring{\(\{4,5,10\}\)}{\{4,5,10\}}}
\label{sec:first}

Fix \(t\ge2\), set \(c=15t+13\) and put \(P=c+19=15(t+2)+2\).
Write a position uniquely as \(n=kP+15y+z\), where
\(k,y\ge0\), \(0\le z<15\), and \(15y+z<P\). Each round \(k\) consists of \(t+2\) rows of
length \(15\) and the two-letter row \(y=t+2\). Use the words
\[
\begin{array}{@{}rl@{\qquad}rl@{}}
 H & \digits{000011112022313} & H' & \digits{000011113022213}\\
 J_1 & \digits{040013113020213} & J & \digits{030014113020213}\\
 K_0 & \digits{340020113120223} & K_* & \digits{130020113120223}\\
 K & \digits{130020113130224} & E & \digits{13}.
\end{array}
\]
The row word \(W(k,y)\) is given by the following rules. For
\(0\le y<t\), it is \(H\) when \(y\le t-k\) and \(H'\) otherwise.
For \(y=t\) it is \(H,J_1,J\) at \(k=0,k=1,k\ge2\), respectively.
For \(y=t+1\) it is \(K_0,K_*,K\) at \(k=0,1\le k\le t,k\ge t+1\),
respectively. The row \(y=t+2\) is always \(E\).
In rounds \(2,\ldots,t+1\), rows \(t-1,\ldots,0\) change from
\(H\) to \(H'\), respectively.
Define \(G(n)=W(k,y)[z]\). In the first round \(n<c\) gives
\(G(n)=H[n\bmod15]\), and the base check of
Section~\ref{sec:method} holds for \(H\).

To verify the remaining positions, use
\(\Phi(C,D)=\slc{C}{4}{15}\slc{D}{0}{4}\) and
\(\Psi(C)=\slc{C}{4}{15}E\,\digits{00}\). Also put
\(Q_t=J\) for \(t=2\) and \(Q_t=H\) for \(t\ge3\).
Table~\ref{tab:first} gives all the inputs of \eqref{eq:local}.
The first row uses \(z_0=13\), and all others use \(z_0=0\).

\begin{table}[!htbp]
\centering
\small
\setlength{\tabcolsep}{5pt}
\renewcommand{\arraystretch}{1.08}
\begin{tabular}{@{}llll@{}}
\toprule
Condition & \(U\) & \(B\) & \(F\)\\
\midrule
\rowcase{k=0,\ y=t}{H}{H}{0^{15}}
\rowcase{k=0,\ y=t+1}{H}{K_0}{\slc{H}{2}{15}\,\digits{00}}
\rowcase{k=0,\ y=t+2}{K_0}{E}{\digits{00}}
\midrule
\rowcase{y=0,\ 1\le k\le t,\ (t,k)\ne(2,2)}{K_*E}{H}{\Phi(H,H)}
\rowcase{y=0,\ t=k=2}{K_*E}{H}{\Phi(H,J_1)}
\rowcase{y=0,\ k=t+1}{K_*E}{H'}{\Phi(H',Q_t)}
\rowcase{y=0,\ k\ge t+2}{KE}{H'}{\Phi(H',Q_t)}
\midrule
\rowcase{1\le y\le t-2,\ 1\le k\le t-y,\ (y,k)\ne(t-2,2)}{H}{H}{\Phi(H,H)}
\rowcase{y=t-2\ge1,\ k=2}{H}{H}{\Phi(H,J_1)}
\rowcase{y=t-2\ge1,\ k=3}{H}{H'}{\Phi(H',J)}
\rowcase{y=t-2\ge1,\ k\ge4}{H'}{H'}{\Phi(H',J)}
\rowcase{1\le y\le t-3,\ k=t-y+1}{H}{H'}{\Phi(H',H)}
\rowcase{1\le y\le t-3,\ k\ge t-y+2}{H'}{H'}{\Phi(H',H)}
\midrule
\rowcase{y=t-1,\ k=1}{H}{H}{\Phi(H,K_0)}
\rowcase{y=t-1,\ k=2}{H}{H'}{\Phi(J_1,K_*)}
\rowcase{y=t-1,\ k\ge3}{H'}{H'}{\Phi(J,K_*)}
\rowcase{y=t,\ k=1}{H}{J_1}{\Psi(K_0)}
\rowcase{y=t,\ 2\le k\le t+1}{H'}{J}{\Psi(K_*)}
\rowcase{y=t,\ k\ge t+2}{H'}{J}{\Psi(K)}
\rowcase{y=t+1,\ k=1}{J_1}{K_*}{\slc{H}{2}{15}\,\digits{00}}
\rowcase{y=t+1,\ 2\le k\le t}{J}{K_*}{\slc{H}{2}{15}\,\digits{00}}
\rowcase{y=t+1,\ k\ge t+1}{J}{K}{\slc{H'}{2}{15}\,\digits{00}}
\rowcase{y=t+2,\ 1\le k\le t}{K_*}{E}{\digits{00}}
\rowcase{y=t+2,\ k\ge t+1}{K}{E}{\digits{00}}
\bottomrule
\end{tabular}
\caption{Local mex data for \(S=\{4,5,10\}\). Expanding the two choices
of \(Q_t\) gives 26 checks.}\label{tab:first}
\end{table}

For \(y\ge1\), the preceding ten letters come from \(W(k,y-1)\).
At \(y=0\) and \(k\ge1\), they are the suffix of
\(W(k-1,t+1)E\). Since \(K_0\) and \(K_*\) agree from column \(2\)
onwards, this suffix is \(\tail(K_*E)\) for \(1\le k\le t+1\),
and \(\tail(KE)\) thereafter. The large option satisfies
\(n-c=(k-1)P+15y+z+19\).
For \(y\le t-1\) and \(k\ge1\) this gives
\(\Phi(W(k-1,y+1),W(k-1,y+2))\). At \(y=t\), the last
four letters of \(F\) come from the short row \(E\) and the first
two letters of row \(0\) of round \(k\), giving
\(\Psi(W(k-1,t+1))\) for \(k\ge1\). For \(k=0\) only
\(z=13,14\) use this move, and both read a zero in row \(0\).
For \(y=t+1\) the option starts at column \(2\) of row \(0\)
of round \(k\) and ends with the first two letters of row \(1\).
For \(y=t+2\) it reads columns \(2,3\) of row \(1\).
The prefixes of length four of \(H,H'\) are both \(\digits{0000}\),
and those of \(K,K_*\) are both \(\digits{1300}\). These formulas
give the boundary rows of Table~\ref{tab:first}.

In the interior \(1\le y\le t-3\), the row \(B\) and the row
\(W(k-1,y+1)\) both switch from \(H\) to \(H'\) precisely when
\(k>t-y\). The preceding row switches one round later, and the
prefix of \(W(k-1,y+2)\) is constant. At \(y=t-2\), the same
large option instead reads the prefix of \(W(k-1,t)\), which
passes through \(H,J_1,J\). At \(y=0\) this latter situation
occurs exactly when \(t=2\). These alternatives exhaust the
rows, including the empty interior ranges when \(t=2\).
Equation~\eqref{eq:local} and the base check now give \(G=g_T\).

For \(k\ge t+1\) each round is \((H')^t JKE\).
The last two letters of round \(t\) are already \(E\), so \(P\)
is a period of \(g_T\) from \(\lambda=(t+1)P-2\).
At \(\lambda-1\) and \(\lambda-1+P\) the letters are
\(K_*[14]=3\) and \(K[14]=4\). Lemma~\ref{lem:min} gives
\(\ell(g_T)=\lambda=(c^2+21c+8)/15\).

For outcomes, \(H,H'\) agree in their zero positions, as do \(J,J_1\)
and \(K_0,K_*,K\) within their respective groups. The outcome
rounds are therefore identical for \(k\ge1\). Comparing rounds \(0\)
and \(1\), the only zero-status changes are in columns \(1\) and \(11\)
of row \(t\), the last at \(c-2\). Lemma~\ref{lem:min} gives
\(\ell(w_T)=c-1\).
There are \(a=5t+9\) ones in an outcome round for \(k\ge1\).
The identity \((t+2)P-(3t+7)a=1\) gives \(\gcd(P,a)=1\),
so Lemma~\ref{lem:min} gives \(p(w_T)=P\).
Since \(P\) is an eventual period of \(g_T\) and every SG period
is an outcome period, \(P=p(w_T)\le p(g_T)\le P\).

The restriction \(t\ge2\) is necessary: for \(c=28\), direct
recurrence gives \(p(g_T)=3\) and \(\ell(g_T)=118\).

\section{The base set \texorpdfstring{\(\{3,7,10\}\)}{\{3,7,10\}}}
\label{sec:second}

Fix \(t\ge1\), set \(c=90t+9\), and put
\(L=c+35=90t+44\), \(P=3tL+31\). Write
\[
 n=mP+kL+45y+z,\qquad m\ge0,\qquad 0\le z<45.
\]
Within each block \(m\), rounds \(0\le k<3t\) have length \(L\):
rows \(0\le y<2t\) have length \(45\) and row \(2t\) has
length \(44\). Round \(k=3t\) consists of a single row of length
\(31\), indexed by \(y=0\). In each row, \(z\) ranges from \(0\)
to one less than its length. These restrictions make the coordinates unique.
Set \(\mu=0\) for \(m=0\) and \(\mu=1\) for \(m\ge1\).
In this section let
\[
 H=\digits{000111022133220031100211322203310001112223332},
 \qquad N=\slc{H}{14}{45}\slc{H}{0}{14}.
\]
The base check of Section~\ref{sec:method} holds for this \(H\)
with \(S=\{3,7,10\}\).

For \(k=3j+\rho\), \(0\le\rho<3\), put
\(f(k)=2t-2j-\rho\) and \((\phi_0,\phi_1,\phi_2)=(13,28,43)\).
The word \(R_\rho^\mu\) is
\(\slc{H}{0}{\phi_\rho}\slc{N}{\phi_\rho}{45}\) with the
following letter replacements:
\[
\begin{array}{@{}ccll@{}}
\toprule
\rho & \phi_\rho & \mu=0 & \mu=1\\
\midrule
0 & 13 & \varnothing & 8\mapsto4,\ 11\mapsto2\\
1 & 28 & 30\mapsto3 & 21\mapsto4,\ 24\mapsto2,\ 27\mapsto3\\
2 & 43 & 43\mapsto3 & 40\mapsto4\\
\bottomrule
\end{array}
\]
For a word \(V\) of length at least \(44\), define \(\mathcal J(V)\) to
be \(\slc{V}{0}{44}\) with columns \(30,37,40\) replaced by
\(4,2,3\), respectively. Write \(J=\mathcal J(N)\),
\(J_\mu=\mathcal J(R_0^\mu)\), and \(D_\mu=\slc{R_0^\mu}{0}{31}\).

Define \(G(n)\) from the following row words. In a full round,
a row \(y<2t\) is \(H\), \(R_\rho^\mu\), or \(N\) according as
\(y<f(k)\), \(y=f(k)\), or \(y>f(k)\). Row \(2t\) is \(J_\mu\)
at \(k=0\) and \(J\) at \(k\ge1\). The short round \(k=3t\)
has word \(D_\mu\).
The splice position in round \(k\) is
\(45f(k)+\phi_\rho=90t+13-30k\); within a block it moves
\(30\) positions towards the start per round.
For \(n<c\) we are in round \(0\) of block \(0\), whose rows
before row \(2t\) are \(H\). In row \(2t\), the columns
\(0,\ldots,8\) of \(J_0\) agree with \(H\), so
\(G(n)=H[n\bmod45]\) for all \(n<c\).

Write \(\Phi(C,D)=\slc{C}{35}{45}\slc{D}{0}{35}\) in this section.
Table~\ref{tab:second} gives \(U,B,F\) for all rows meeting \(n\ge c\).
The first row uses \(z_0=9\) and all others use \(z_0=0\).
In the table \(f=f(k)\), and for \(m\ge1\) let \(\nu=0\) when
\(m=1\) and \(\nu=1\) when \(m\ge2\).

\begin{table}[!htbp]
\centering
\small
\setlength{\tabcolsep}{5pt}
\renewcommand{\arraystretch}{1.12}
\begin{tabular}{@{}llll@{}}
\toprule
Condition & \(U\) & \(B\) & \(F\)\\
\midrule
\rowcase{k=m=0,\ y=2t}{H}{J_0}{0^9\slc{H}{0}{35}}
\rowcase{k=0,\ m\ge1,\ y\le2t-2}{H}{H}{\Phi(H,H)}
\rowcase{k=0,\ m\ge1,\ y=2t-1}{H}{H}{\slc{J}{21}{44}\slc{D_\nu}{0}{22}}
\rowcase{k=0,\ m\ge1,\ y=2t}{H}{J_1}{\slc{D_0}{22}{31}\slc{H}{0}{35}}
\rowcase{k=3t}{J}{D_\mu}{\slc{R_2^\mu}{35}{45}\slc{N}{0}{21}}
\midrule
\multicolumn{4}{@{}l@{}}{\emph{For the following rows, \(1\le k<3t\).}}\\
\rowcase{y<f}{H}{H}{\Phi(H,H)}
\rowcase{y=f,\ \rho=0}{H}{R_0^\mu}{\Phi(R_2^\mu,N)}
\rowcase{y=f,\ k=1}{H}{R_1^\mu}{\Phi(H,J_\mu)}
\rowcase{y=f,\ \rho=1,\ k\ge4}{H}{R_1^\mu}{\Phi(H,R_0^\mu)}
\rowcase{y=f,\ \rho=2}{H}{R_2^\mu}{\Phi(H,R_1^\mu)}
\rowcase{f<y\le2t-2,\ (y,\rho)\ne(f+1,2)}{N}{N}{\Phi(N,N)}
\rowcase{y=f+1\le2t-2,\ \rho=2}{R_2^\mu}{N}{\Phi(N,N)}
\rowcase{y=2t-1,\ k=2}{R_2^\mu}{N}{\Phi(N,J)}
\rowcase{y=2t-1,\ k\ge3}{N}{N}{\Phi(N,J)}
\rowcase{y=2t}{N}{J}{\slc{J}{35}{44}\slc{H}{0}{35}}
\bottomrule
\end{tabular}
\caption{Local mex data for \(S=\{3,7,10\}\). There are 22 checks
after expanding the choices of \(\mu,\nu\) and counting the
repeated \((H,H,\Phi(H,H))\) entry once.}\label{tab:second}
\end{table}

To check coverage, first consider the small options. When \(y>0\)
the preceding ten letters are the suffix of the previous row.
At \(y=0\) they come from the previous round or from the previous
block's short round. The equalities
\[
\begin{split}
\tail(J)=\tail(J_\mu)=\tail(D_\mu)&=\tail(H),\\
\tail(R_0^\mu)=\tail(R_1^\mu)&=\tail(N),\qquad
\slc{R_2^\mu}{0}{35}=\slc{H}{0}{35}
\end{split}
\]
follow from the words, so they account for the choices of \(U\).

For the large option write \(X=45y+z\).
If \(k\ge1\), \(n-c\) is at offset \(X+35\) from the start of
round \(k-1\) in block \(m\).
For a full row \(y<2t\) in a full round this gives
\(\Phi\) of rows \(y,y+1\) in round \(k-1\).
When \(y<f\), both are old-phase rows except possibly
row \(y+1=f(k-1)\). In that case \(\rho=0\) and the
prefix of \(R_2^\mu\) agrees with \(H\).
When \(y=f\), the previous front is at \(f(k-1)=f\) if
\(\rho=0\), or \(f+1\) if \(\rho\ne0\). The two rows read
are consequently \(R_2^\mu,N\) or \(H,R_{\rho-1}^\mu\).
The latter has \(J_\mu\) in place of \(R_{\rho-1}^\mu\)
when \(k=1\), since the previous front is then in row \(2t\).
When \(y>f\), the suffix read is that of \(N\), or, for
\(\rho\ne0\) and \(y=f+1\), of \(R_{\rho-1}^\mu\) with the
same suffix. The second word
read is \(N\) except at \(y=2t-1\), where it is \(J\).
In row \(2t\), offset \(X+35\) reads columns \(35,\ldots,43\)
of the previous round's last row and columns \(0,\ldots,34\)
of the current round's first row. The suffix
\(\slc{J_\mu}{35}{44}\) agrees with \(\slc{J}{35}{44}\),
and the first row has prefix \(\slc{H}{0}{35}\), also when it is
\(R_2^\mu\) at \(k=3t-1\).
In the short round \(X+35\) reads the suffix of \(R_2^\mu\)
and the first \(21\) letters of \(N\) from the last full round.

For \(k=0,m\ge1\), use \(P-c=(3t-1)L+66\).
The large option starts at offset \(45(y+1)+21\) from the start
of the previous block's last full round, then advances \(z\) places.
For \(y\le2t-2\) the word read is
\(\slc{N}{21}{45}\slc{N}{0}{21}=\Phi(H,H)\), using
\(\slc{J}{0}{21}=\slc{N}{0}{21}\) when \(y+2=2t\).
For \(y=2t-1\) the option crosses into the previous short
round, and for \(y=2t\) it crosses from that short round into
the current block. Here \(\slc{D_0}{22}{31}=\slc{D_1}{22}{31}\).
In the remaining case \(k=m=0\), only row \(2t\) at \(z\ge9\)
reads \(c\)-options, at \(n-c=z-9\). This proves coverage of
Table~\ref{tab:second} also for \(t=1\), where the intermediate
row ranges that are empty contribute no conditions. By
\eqref{eq:local}, \(G=g_T\).

The blocks \(m\ge1\) are identical. In the short round of
block \(0\), \(D_0\) and \(D_1\) differ only at columns \(8\)
and \(11\), with \(D_0[11]=3\) and \(D_1[11]=2\).
Thus \(P\) is a period from \(\lambda=3tL+12=P-19\),
and \(\lambda-1\) is a mismatch. Lemma~\ref{lem:min} gives
\(\ell(g_T)=P-19\).

Each change between \(R_\rho^0\) and \(R_\rho^1\) preserves zero
status, as do the changes between \(J_0,J_1\) and \(D_0,D_1\).
Hence \(w_T\) repeats from \(n=0\) with period \(P\).
To prove minimality, count zeros and adjacent zero pairs in the
SG words. The words \(H,N,R_\rho^\mu\) each have \(12\) zeros
and \(6\) adjacent zero pairs. Each of \(J,J_\mu\) has \(11\)
zeros and \(6\) such pairs, and each \(D_\mu\) has \(8\) zeros
and \(4\) pairs. All these words end in \(2\), so no zero pair
crosses a row boundary, including the wrap at the end of a block.
The counts of the cyclic factors \(1\) and \(11\) in an
outcome block are therefore
\[
 a=3t(24t+11)+8,\qquad b=3t(12t+6)+4.
\]
The identity \((6t+4)a-12tb-P=1\) gives \(\gcd(P,a,b)=1\).
Lemma~\ref{lem:min} gives \(p(w_T)=P\), which forces \(p(g_T)=P\).
Finally, substitution of \(c=90t+9\) gives
\(30P=c^2+26c+615\), completing the proof of
Theorem~\ref{thm:main}.

\section*{Data and code availability}
All word data and the finite mex test are included in the article.
The supplementary script \texttt{verify\_quadratic\_subtraction\_note.py}
checks all 48 local identities, compares the constructions with direct
mex recurrence, and independently computes the outcome sequence.
Repeated states, primitive cycle words, and boundary mismatches certify
the periods and preperiods of each tested instance.

\section*{Declaration of generative AI and AI-assisted technologies}
ChatGPT was used to assist with language polishing, LaTeX preparation,
and calculations. The authors are responsible for the manuscript and
its mathematical claims.

\end{document}